\documentclass[oneside,10pt]{article}          
\usepackage[b5paper]{geometry}	               
\usepackage{amsfonts,amsmath,latexsym,amssymb} 
\usepackage{theorem}                           
\usepackage{mathrsfs,upref}                    
\usepackage{mathptmx}		                   
\usepackage{arXiv}	                           

\newtheorem{theorem}{Theorem}
\newtheorem{statement}{Statement}
\newtheorem{proposition}{Proposition}
\newtheorem{lemma}{Lemma}

\theoremstyle{definition}

\newtheorem{example}{Example}

\newtheorem{conjecture}{Conjecture}

\usepackage[dvipdf]{graphicx}

\newfont{\micros}{msam10 scaled 700}

\begin{document}

\title[On the Erd\" os-Mordell Inequality for Triangles in the Taxicab Geometry]{On the Erd\" os-Mordell Inequality for Triangles in Taxicab Geometry}


\author{Maja Petrovi\' c, Branko Male\v sevi\' c and Bojan Banjac}

\address{Maja Petrovi\' c, University of Belgrade - The Faculty of Transport and Traffic Engineering, \\ Vojvode Stepe $305$, Belgrade, Serbia \\
\email{majapet@sf.bg.ac.rs}}

\address{Branko Male\v sevi\' c, University of Belgrade - School of Electrical Engineering, \\ Bulevar Kralja Aleksandra $73$, Belgrade, Serbia \\
\email{malesevic@etf.bg.ac.rs}}

\address{Bojan Banjac, University of Novi Sad - Faculty of Technical Sciences, \\  Trg Dositeja Obradovi\' ca $6$, Novi Sad, Serbia \\
\email{bojan.banjac@uns.ac.rs}}

\CorrespondingAuthor{Branko Male\v sevi\' c}


\date{03.12.2018}                               




\begin{abstract}
In this work the {\sc Erd\" os-Mordell}'s inequality is examined for the case of a triangle $ABC$ in the taxicab plane geometry.
It is shown that the {\sc Erd\" os-Mordell}'s~\mbox{inequality} \mbox{$R_A + R_B + R_C \geq w (r_a + r_b + r_c)$} holds for triangles
with appropriate positions for its points \mbox{$A$, $B$} and $C$, if $w = 3/2$.
\end{abstract}

\maketitle

\section{Introduction}

Let the distance between two points, as well as the distance between a line and a point be defined in the Euclidean plane.
Then, for a triangle $ABC$ in such a plane the {\sc Erd\" os-Mordell}'s inequality holds \cite{Erdos_1935}, \cite{Mordell_1937}:
\begin{equation}
    \label{E1}
    R_A + R_B + R_C
    \geq
    2\left(r_a + r_b + r_c\right)
    \smallskip
\end{equation}
where $R_A$, $R_B$ and $R_C$ are distances from the interior point $M$ of $\triangle ABC$ to vertices $A$, $B$ and $C$ respectively
and $r_a$, $r_b$ and $r_c$ are distances from the point $M$ of the triangle to the corresponding edges which contain the vertices of $\triangle\,ABC$ (Fig. 1).

\begin{center} 

\vspace*{2.0mm} \hspace*{-7.0mm} \includegraphics*[height=45.0mm]{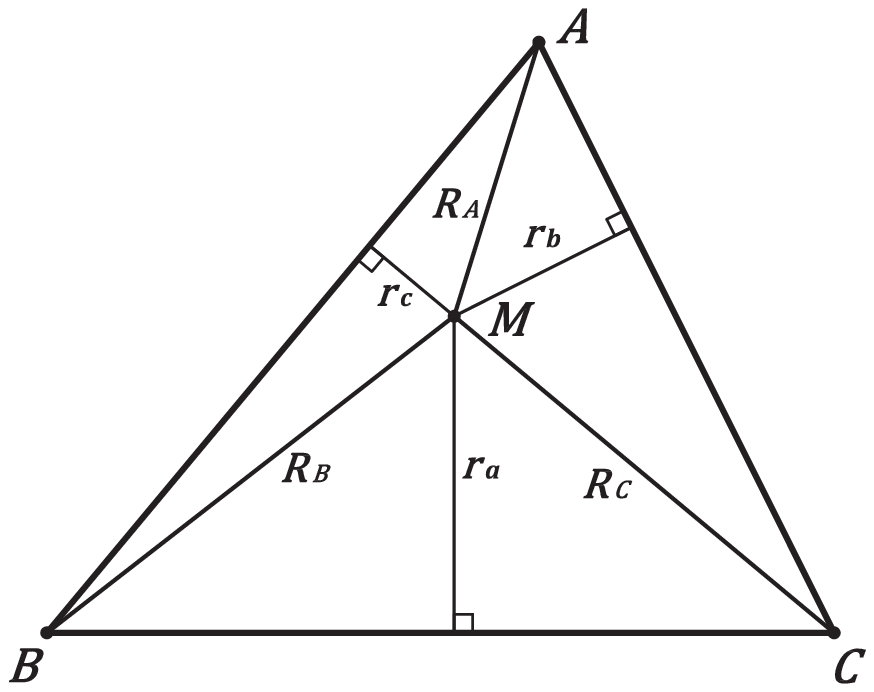}

\smallskip
\noindent
\textit{Figure 1: A geometric illustration of the {\sc Erd\" os}-{\sc Mordell} inequality in $\triangle ABC$}
\end{center}

\break
\noindent
Let there be two points, $A\left(x_A,y_A\right)$ and $B\left(x_B,y_B\right)$, then the distance between them in taxicab geometry is defined as:
\begin{equation}
d_1\left(A,B\right) = |x_A-x_B|+|y_A-y_B|.
\end{equation}
This distance is also called the Manhattan or city block distance. This metric is a special case of the Minkowski metric of order $k$ (where $k\geq 1$)
which is defined by the following formula:
\begin{equation}
    \label{E2}
    d_k\left(A,B\right) = \left(|x_A-x_B|^k+|y_A-y_B|^k \right)^ \frac{1}{k}
    \smallskip
\end{equation}
The Minkowski metric contains in itself the taxicab metric for the value $k = 1$ and the Euclidean metric for $k = 2$ \cite{Kraus_1986}.
The term "taxicab" was first introduced by {\sc K. Menger} \cite{Menger_1952}. A graphical representation of distances between points
$A$ and $B$ is given in Fig.~2, in taxicab metric with $d_1$ (dashed/long dashed lines) and in Euclidean metric with $d_2$ (continuous line).

\begin{center} 

\vspace*{0.0mm} \hspace*{-8.0mm} \includegraphics*[height=50.0mm]{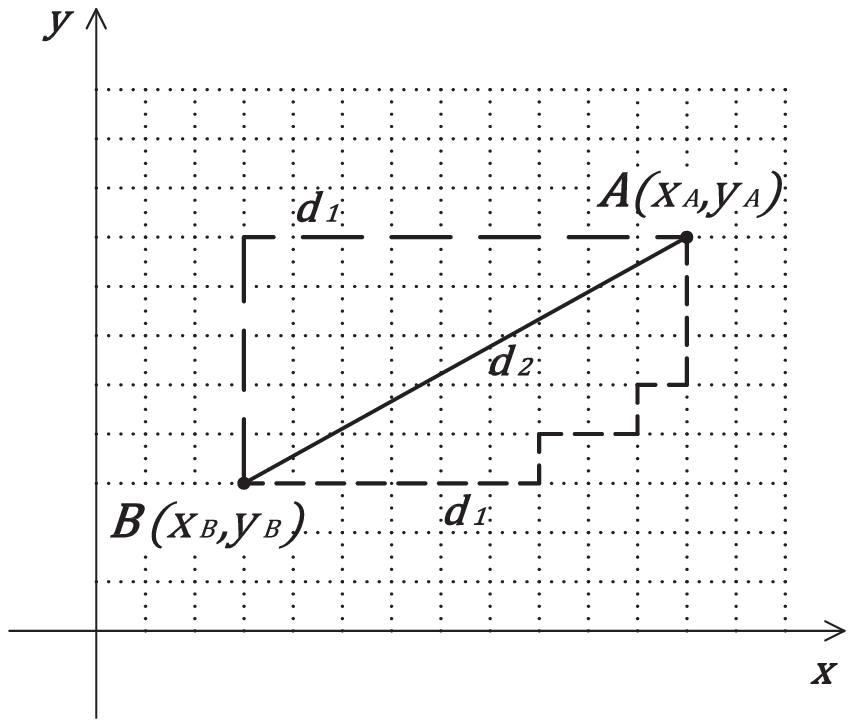}

\noindent
\textit{Figure 2: A geometric illustration of the Minkowski and the Euclidean distances \\ between two points}
\end{center}

In the rest of this work, only taxicab distances are considered.

\smallskip
Let the $\triangle ABC$ be a triangle with vertices $A\left(0,r\right), \, B\left(p,0\right), \, C\left(q,0\right), \, p \ne q, \, r\ne 0$.
Without diminishing generality, let $p<q$. We denote by $M\left(x,y\right)$ an arbitrary point in the plane of the triangle $\triangle ABC$  (Fig. 1).
The Taxicab distance from the point $M$ to the points $A$, $B$ and $C$, are given by functions:
\begin{equation}
\begin{split}
\label{E3}
     R_A = d_1\left(M,A\right) = |x|+|y-r|,
    \\
    R_B = d_1\left(M,B\right) = |x-p|+|y|,
    \\
    R_C = d_1\left(M,C\right) = |x-q|+|y|.
\end{split}
\end{equation}

Recently, general formulae for distance in taxicab geometry were analyzed in the paper \cite{Colakoglu_2019}.
Authors {\sc Kaya} et al. \cite{Kaya_at_all_2000} define the distance of a point to a line in taxicab plane geometry with the following statement:
\begin{lemma}
Distance of point $M\left(x_M,y_M\right)$ to the line $\ell\!:ax+by+c=0$ in the Taxicab plane is$:$
\begin{equation}
    \label{E4}
    d_1\left(M,\ell \right)=
    \frac{|ax_M+by_M+c|}{max \{ |a|,|b| \} }.
\end{equation}
\end{lemma}
Let us notice that
\begin{equation}
\label{E4_ell}
r_{a} = d_1(M,\ell_{\!\mbox{\tiny $BC$}}),\;
r_{b} = d_1(M,\ell_{\!\mbox{\tiny $AC$}}),\;
r_{c} = d_1(M,\ell_{\!\mbox{\tiny $AB$}}).
\end{equation}
Based on (\ref{E3}) and (\ref{E4_ell}), the {\sc Erd\" os-Mordell}'s inequality (\ref{E1}) for $\triangle\,ABC$ in taxicab metric is defined by the following relation:
\begin{equation}
    \label{E5}
    |x|+|y-r|+|x-p|+|x-q|+2|y|\ge
    2\left(|y| + \frac{|qr-rx-qy|}{max \{|r|,|q|\}}+\frac{|pr-rx-py|}{max \{|r|,|p| \}}\right).
\end{equation}
Inequalities in the taxicab geometry are the topic of recent research, see e.g \cite{Kocayusufoglu_2006}.
Let us emphasize that the topic of the {\sc Erd\" os-Mordell} inequality is current, as it has been shown in the papers \cite{Dergiades_2004},
\cite{Ghandehari_2018}, \cite{Liu_2015} --  \cite{Liu_2019}, \cite{Walker_2016}  and books \cite{Bottema_at_all_1969} and \cite{Mitrinovic_at_all_1989}.
\mbox{\sc V.~Pambuccian} proved that, in the plane of absolute geometry, the {\sc Erd\" os-Mordell} inequality
is an equivalent to the non-positive curvature \cite{Pambuccian_2008}. In the paper \cite{Malesevic_at_all_2014} is given an extension of
the {\sc Erd\" os-Mordell} inequality on the interior of the {\sc Erd\"os-Mordell} curve. In relation to
the {\sc Erd\" os-Mordell} inequality {\sc N.~Dergiades} in the paper \cite{Dergiades_2004} proved one extension of
the {\sc Erd\" os-Mordell} type inequality. Most notably, the {\sc Erd\" os-Mordell} inequality has been considered in the taxicab plane geometry by {\sc N. S\" onmez}
who has shown that (\ref{E1} is a strict inequality: $R_A + R_B + R_C > 2\left(r_a + r_b + r_c\right)$, \cite{Sonmez_2009}. In this work we prove that the conclusion reached by
{\sc N. S\" onmez} is incorrect. That shall be shown through the following example.
\begin{example} (counterexample)
Let the vertices of $\triangle ABC$ be given with $p=-20,\, q=40,\, r=30$ and let point $M(0,m)$ be defined with $m = 2$ (Fig. 3).
The taxicab distance from the point $M$ to the vertices of $\triangle ABC$ is given by (\ref{E3}) and the distance from point $M$
to the lines $\ell_{\!\mbox{\tiny $BC$}}\!:y=0,\, \ell_{\!\mbox{\tiny $AC$}}\!:-rx-qy+qr=0$ and $\ell_{\!\mbox{\tiny $AB$}}\!:-rx-py+pr=0$ is given by (\ref{E4}):
{\small \begin{equation}
   \begin{split}
   \label{E6}
    R_A &= d_1\left(M,A\right) = 28, \,\,\, \,\,\, \,\,\,
    R_B = d_1\left(M,B\right) = 22, \,\,\, \,\,\,
    R_C = d_1\left(M,C\right) = 42,                                                \\[-0.5 ex]
    r_a &= d_1\left(M,\ell_{\!\mbox{\tiny $BC$}} \right)= 2, \,\,\, \,\,\,   \,\,\, \,\,\,\,
    r_b = d_1\left(M,\ell_{\!\mbox{\tiny $AC$}} \right)=  28, \,\,\, \,\,\,   \,\,\,
    r_c = d_1\left(M,\ell_{\!\mbox{\tiny $AB$}} \right)=  \mbox{\small $\displaystyle\frac{56}{3}$}.
   \end{split}
 \end{equation}}

\begin{center} 

\vspace*{-2.0mm} \hspace*{0.0mm} \includegraphics*[height=45.0mm]{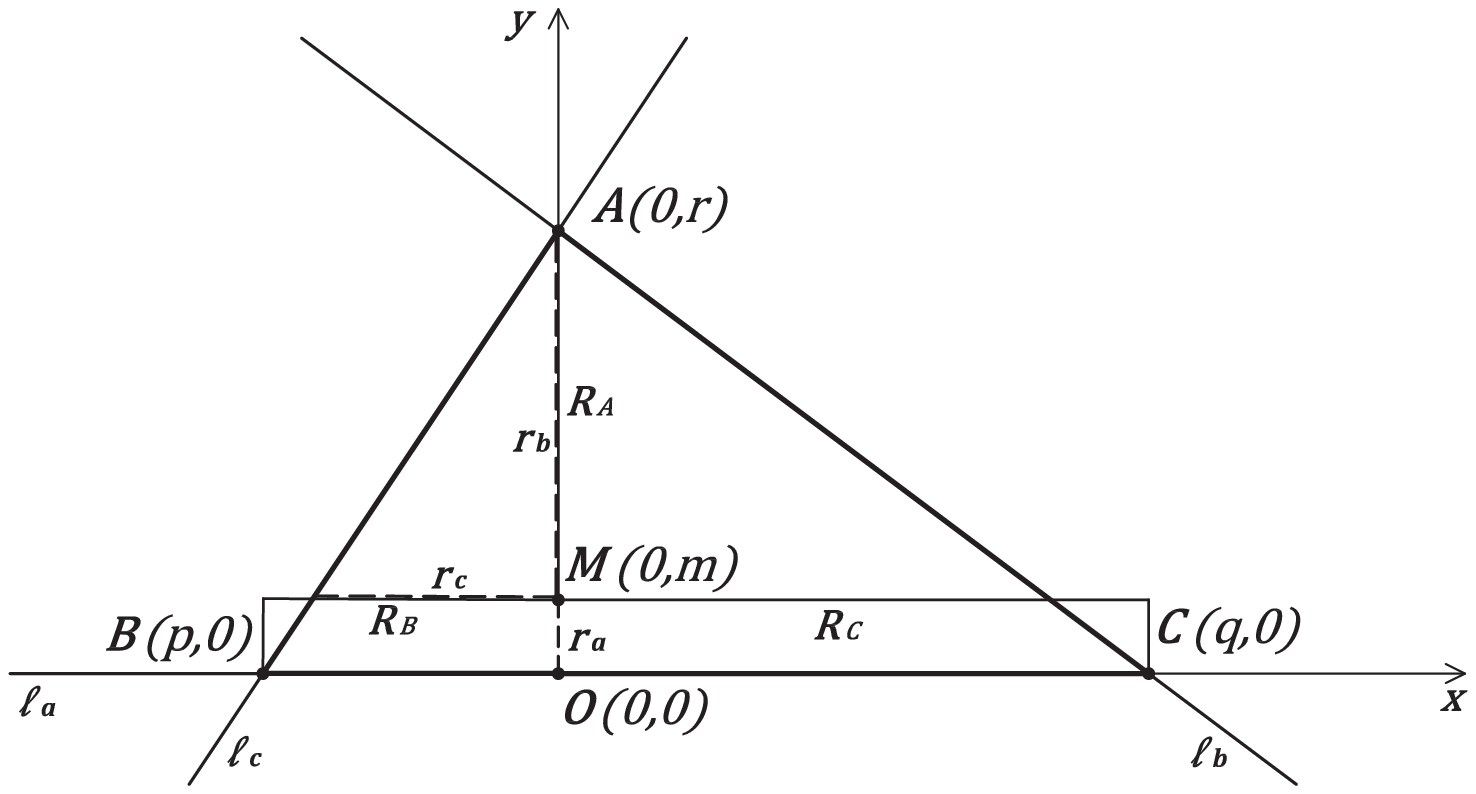}

\noindent
\textit{Figure 3: A geometric illustration of the counterexample}
\end{center}

\break

\bigskip
From (\ref{E6}) we obtain $L=R_A + R_B + R_C=92$ and $R=r_a + r_b + r_c=\frac{146}{3}$. In the case of the {\sc Erd\" os-Mordell} inequality, it holds that
$L \geq 2R$ i.e $92 \geq 97.\overline{3}\,$. From this follows that the {\sc Erd\" os-Mordell} inequality does not hold for all interior points of $\triangle\,ABC$. \hfill $\Box$
\end{example}

In the rest of this paper, the {\sc Erd\" os-Mordell} inequality in taxicab geometry is considered in the form:
\begin{equation}
    \label{E7}
    R_A + R_B + R_C \ge
    w\left(r_a + r_b + r_c\right),
    \smallskip
\end{equation}
where the positive real number w is defined as such that the previous inequality holds for all interior points of $\triangle ABC$.$\,$The main goal of the work is to,
for all positive values of the weight coefficient $w$, determine a upper bound $\mathfrak{M}$ such that the {\sc Erd\" os-Mordell} inequality holds for $0 < w \leq \mathfrak{M}$.

\section{The Main Results}

The {\sc Erd\" os-Mordell} inequality in taxicab plane geometry has the following form:
\begin{equation}
    \label{E8}
\!|x|+|y-r|+|x-p|+|x-q|+2|y|\ge
    w\!\left(|y| + \displaystyle\frac{|qr-rx-qy|}{max\{|r|,|q|\}} + \displaystyle\frac{|pr-rx-py|}{max\{|r|,|p|\}}\right)\!.
    \smallskip
\end{equation}

It should be noted that the {\sc Erd\" os-Mordell} inequality in the taxicab plane geometry defined by (\ref{E8}) refers to triangles ABC with the appropriate positions of points
$A(0,r)$, $B(p,0)$ and $C(q,0)$ in two cases. The first case is when coordinates $p$, $q$ and $r$ are positive and the second case is when the $p$ coordinate is negative, with positive
$q$ and $r$ coordinates. Furthermore, we do not consider the general position of the triangle in the taxicab plane nor the rotation of such a triangle to $\triangle ABC$.

\smallskip
{\bf $1^{\circ}$}  We analyze $\triangle\,ABC$ with $p, q, r > 0$ (see Fig. 4), then, for all interior points of the triangle holds:
\begin{equation}
\begin{split}
    \label{E9}
    |x|=x,  \,\,\, \,\,\, \,\,\, \,\,\, \,\,\,
    |x-p|= \left\{
           \begin{array}{ccc}
                    p-x      & : &   x<p           \\[1.0 ex]
                    x-p      & : & x\geq p
           \end{array} ,
           \right. \,\,\, \,\,\, \,\,\, \,\,\, \,\,\,
    |x-q|=q-x,
\\
    |y|=y, \,\,\, \,\,\, \,\,\, \,\,\, \,\,\,
    |y-r|=r-y,
\\
    |qr-rx-qy|=qr-rx-qy, \,\,\, \,\,\, \,\,\, \,\,\, \,\,\,
    |pr-rx-py|=-pr+rx+py.
    \smallskip
\end{split}
\end{equation}

Then, the form of the {\sc Erd\" os-Mordell} inequality (\ref{E8}) becomes:
\begin{equation}
    \label{E10}
    \left\{
           \begin{array}{ccc}
                    q+r+y+p-x \ge
    w\left(y + \displaystyle\frac{qr-rx-qy}{max\{r,q\}} + \displaystyle\frac{-pr+rx+py}{max\{r,p\}}\right)      & : &   x<p           \\[2.5 ex]
                    q+r+y+x-p  \ge
    w\left(y + \displaystyle\frac{qr-rx-qy}{max\{r,q\}} + \displaystyle\frac{-pr+rx+py}{max\{r,p\}}\right)     & : & x\geq p
           \end{array}
           \right.
    \smallskip
\end{equation}
Symmetric positions of $\triangle ABC$ relative to the coordinate axes can be analogously considered.

\medskip
{\bf $2^{\circ}$}
We analyze  $\triangle\,ABC$ with $p<0$ \, and  \, $q,r>0$ (see Fig. 4), then, for all interior points of the triangle holds:
\begin{equation}
\begin{split}
    \label{E11}
    |x|=\left\{
           \begin{array}{ccc}
                    -x      & : &   x<0           \\[1.0 ex]
                    x      & : & x\geq 0
           \end{array} ,
           \right. \,\,\, \,\,\, \,\,\, \,\,\, \,\,\,
    |x-p|= x-p,  \,\,\, \,\,\, \,\,\, \,\,\, \,\,\,
    |x-q|=q-x,
\\
    |y|=y, \,\,\, \,\,\, \,\,\, \,\,\, \,\,\,
    |y-r|=r-y,
\\
    |qr-rx-qy|=qr-rx-qy,  \,\,\, \,\,\, \,\,\, \,\,\, \,\,\,
    |pr-rx-py|=-pr+rx+py.
    \smallskip
\end{split}
\end{equation}

Then, the form of the {\sc Erd\" os-Mordell} inequality (\ref{E8}) becomes:
\begin{equation}
    \label{E12}
    \left\{
           \begin{array}{ccc}
                    -p+q+r+y-x\ge
    w\left(y + \displaystyle\frac{qr-rx-qy}{max\{r,q\}} + \displaystyle\frac{-pr+rx+py}{max\{r,-p\}}\right)      & : &   x<0           \\[2.5 ex]
                    -p+q+r+y+x\ge
    w\left(y + \displaystyle\frac{qr-rx-qy}{max\{r,q\}} + \displaystyle\frac{-pr+rx+py}{max\{r,-p\}}\right)      & : & x\geq 0
           \end{array}
           \right.
    \smallskip
\end{equation}
As in case {\bf $1^{\circ}\!\!$}, symmetric positions of $\triangle ABC$ relative to the coordinate axes can be analogously considered.

\medskip
Let us notice that for point $A(0,r)$, there exist the following subcases:
{\small $$
\!\!\!\!
\begin{array}{clll}
1^{\circ}
&
\langle \textbf{\textit{a}} \rangle  \;  0 < r\leq p < q,
&
\langle \textbf{\textit{b}} \rangle \;  0\leq p < r \leq q,
&
\langle \textbf{\textit{c}} \rangle \;  0 \leq p < q < r;                                         \\[2.0 ex]
\begin{array}{l}
\textrm{\footnotesize For this subcase, see Fig. $4/1^{\circ}$} \\[-0.5 ex]
\textrm{\footnotesize with representations:}
\end{array}
&
\langle \textbf{\textit{a}} \rangle \; \mbox{$\begin{array}{l}
                                              \textrm{\footnotesize long and double-short} \\[-0.5 ex]
                                              \textrm{\footnotesize dashed line,}
                                              \end{array}$}
&
\langle \textbf{\textit{b}} \rangle \; \textrm{\footnotesize dashed line,}
&
\langle \textbf{\textit{c}} \rangle \; \textrm{\footnotesize continuous line};                    \\[4.0 ex]
2^{\circ}
&
\langle \textbf{\textit{a}} \rangle  \;  0 < r\leq -p \leq q,
&
\langle \textbf{\textit{b}} \rangle \;  0 < -p \leq r \leq q,
&
\langle \textbf{\textit{c}} \rangle \;  0 < -p \leq q < r;                                        \\[2.0 ex]
\begin{array}{l}
\textrm{\footnotesize For this subcase, see Fig. $4/2^{\circ}$} \\[-0.5 ex]
\textrm{\footnotesize with representations:}
\end{array}
&
\langle \textbf{\textit{a}} \rangle \; \mbox{$\begin{array}{l}
                                              \textrm{\footnotesize long and double-short} \\[-0.5 ex]
                                              \textrm{\footnotesize dashed line,}
                                              \end{array}$}
&
\langle \textbf{\textit{b}} \rangle \; \textrm{\footnotesize dashed line,}
&
\langle \textbf{\textit{c}} \rangle \; \textrm{\footnotesize continuous line}.
\end{array}
$$}

\begin{center} 

\vspace*{2.0mm} \hspace*{-2.0mm} \includegraphics*[height=57.5mm]{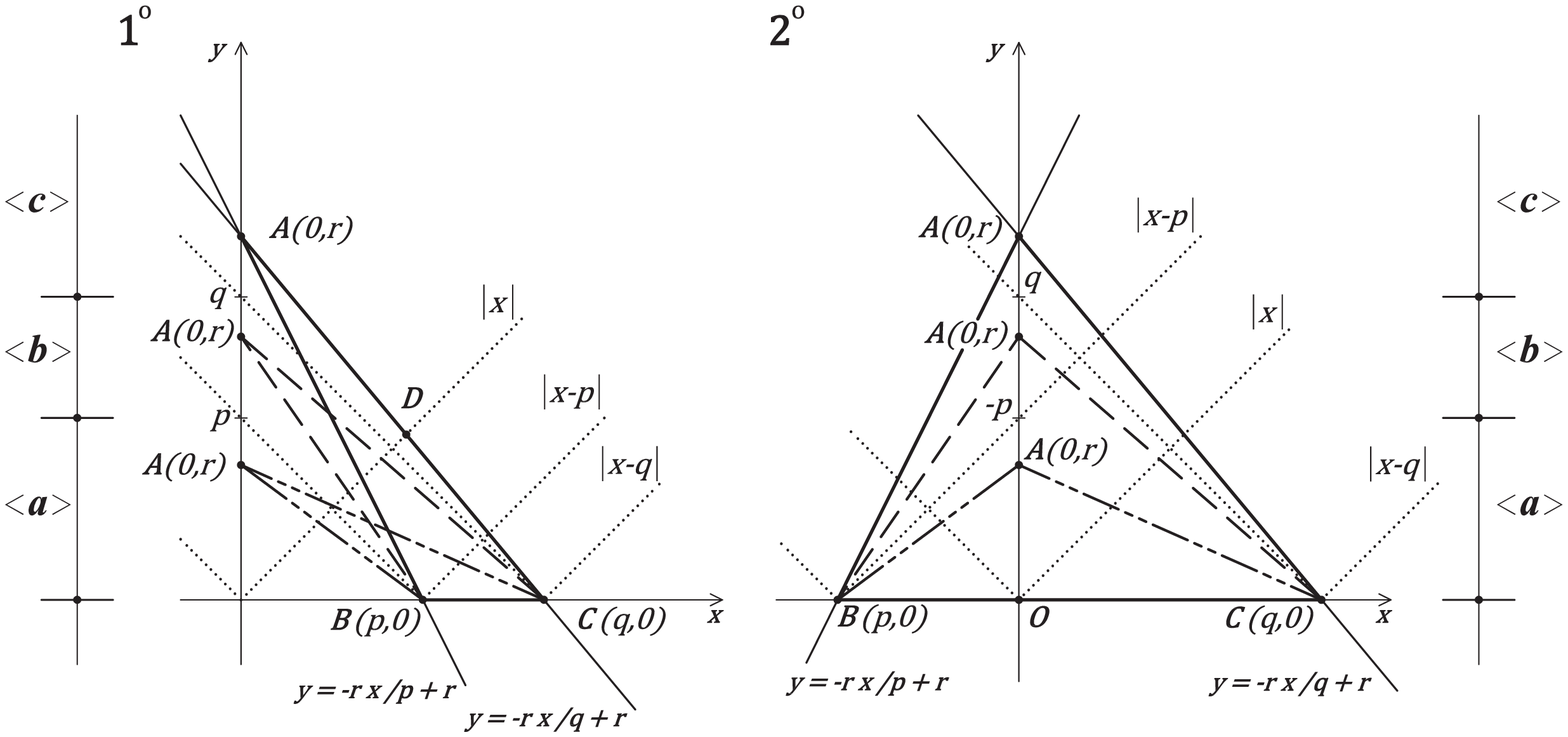}

\smallskip
\noindent
\textit{Figure 4: The two types of triangles $ABC$ with subcasses}
\end{center}

\smallskip
In formula (\ref{E9}), for the first triangle type $(i = 1)$, branching is achieved for $x = p$, where $p$ will then be denoted with $x_1.\,$In
formula (\ref{E11}), for the second triangle type $(i = 2)$, branching is achieved for $x = 0$, where $0$ will then be denoted with $x_2$.
Then, the {\sc Erd\" os-Mordell} inequality (\ref{E8}), with weight coefficient $w > 0$, is considered with the following theorem:

\begin{theorem}{\label{theor}}
 It holds:
\begin{equation}
    \!\label{E13}
    R_A + R_B + R_C \!\ge\!w\left(r_a + r_b + r_c\right)\,
    \Longleftrightarrow
     \left\{
           \begin{array}{cccc}
                    \!\! \alpha_{i1} x + \beta_{i1} y + \gamma_{i1} \ge 0  & : & x<x_i     & \left[\Pi_{i1}\right]      \\[1.0 ex]
                    \!\! \alpha_{i2} x + \beta_{i2} y + \gamma_{i2} \ge 0  & : & x\geq x_i & \left[\Pi_{i2}\right]
           \end{array}
           \right.
       \smallskip
\end{equation}
where coefficients $\alpha_{ij}, \beta_{ij}, \gamma_{ij}$\,$(j\!=\!1,2)$, are given by Tab.$\,$1 for $i \!=\! 1$ and Tab.$\,$2~for~\mbox{$i \!=\! 2$}.

\end{theorem}

\begin{center} 
\begin{tabular}[t]{|c|c|r|r|r|}
\hline
& &\multicolumn{3}{p{7.5cm}|}%
{\centering  $\Pi_{1j} : \alpha_{1j} x + \beta_{1j} y + \gamma_{1j} \ge 0$} \\
\cline{3-5}
\multicolumn{1}{|c|}{$1^{\circ}$}
& \multicolumn{1}{c|}{}
& \multicolumn{1}{c|}{\small  $\langle \textbf{\textit{a}} \rangle$ }
& \multicolumn{1}{c|}{\small  $\langle \textbf{\textit{b}} \rangle$ }
& \multicolumn{1}{c|}{\small  $\langle \textbf{\textit{c}} \rangle$ }\\
\multicolumn{1}{|c|}{}
& \multicolumn{1}{c|}{}
& \multicolumn{1}{c|}{\scriptsize   $0 < r \leq p < q$}
& \multicolumn{1}{c|}{\scriptsize   $0 \leq p < r \leq q$}
& \multicolumn{1}{c|}{\scriptsize   $0 \leq p < q < r$}\\
\hline
\hline
\small $ \!\!$  & \!\!\! \mbox{\small $\,\alpha_{11}$}
                & \footnotesize $(p-q)wr\!-\!pq$ & \footnotesize $r(w(r-q)\!-\!q)$ & \small $-r$
\\
 $\Pi_{11}$ & \!\!\! \mbox{\small $\,\beta_{11}$}
                     & \footnotesize $-pq(w-1)$ & \footnotesize $q(r-pw)$ & \footnotesize $w(q\!-\!p\!-\!r)\!+\!r$
\\
\scriptsize $x<p$  & \!\! \mbox{\small $\,\gamma_{11}$}
                & \footnotesize $pq(p+q+r)$ & \footnotesize $qr(w(p-r)+p+q+r)$ & \footnotesize $r(w(p-q)+p+q+r)$
\\
\hline
\small $\!\! $  & \!\!\! \mbox{\small $\,\alpha_{12}$}
                & \footnotesize $(p-q)wr+pq$ & \footnotesize $r(w(r-q)+q)$ & \small $r$
\\
 $\Pi_{12}$ & \!\!\! \mbox{\small $\,\beta_{12}$}
                    & \footnotesize $-pq(w-1)$ & \footnotesize $q(r-pw)$ & \footnotesize $w(q\!-\!p\!-\!r)\!+\!r$
\\
\scriptsize $x\geq p$ & \!\! \mbox{\small $\,\gamma_{12}$}
                      & \footnotesize $\,\,pq(-p+q+r)$ & \footnotesize $\,qr(w(p-r)-p+q+r)$ & \footnotesize  $\,r(w(p-q)-p+q+r)$
\\
\hline
\end{tabular}
\end{center}
\begin{center}
  \textit{Table 1: The \textsc{Erd\" os-Mordell} inequality in the  taxicab plane geometry for case $1^{\circ}$}
\end{center}
\begin{center} 
\begin{tabular}[t]{|c|c|r|r|r|}
\hline
& &\multicolumn{3}{p{7.5cm}|}%
{\centering  $\Pi_{2j}: \alpha_{2j} x + \beta_{2j} y + \gamma_{2j} \geq 0$}  \\
\cline{3-5}
\multicolumn{1}{|c|}{$2^{\circ}$}
& \multicolumn{1}{c|}{}
& \multicolumn{1}{c|}{\small  $\langle \textbf{\textit{a}} \rangle$ }
& \multicolumn{1}{c|}{\small  $\langle \textbf{\textit{b}} \rangle$ }
& \multicolumn{1}{c|}{\small  $\langle \textbf{\textit{c}} \rangle$ }\\
\multicolumn{1}{|c|}{}
& \multicolumn{1}{c|}{}
& \multicolumn{1}{c|}{\scriptsize   $0 < r \leq -p \leq q$}
& \multicolumn{1}{c|}{\scriptsize   $0 < -p \leq r \leq q$}
& \multicolumn{1}{c|}{\scriptsize   $0 < -p \leq q < r$}\\
\hline
\hline
\small $ \!\!$  & \!\!\! \mbox{\small $\,\alpha_{21}$}
                & \footnotesize $pq\!-\!(p+q)wr$ & \footnotesize $r(w(r-q)\!-\!q)$ & \small $-r$
\\
 $\Pi_{21}$ & \!\!\! \mbox{\small $\,\beta_{21}$}
                     & \footnotesize $-pq(w+1)$ & \footnotesize $q(r-pw)$ & \footnotesize $w(q\!-\!p\!-\!r)\!+\!r$
\\
\scriptsize $x<0$  & \!\! \mbox{\small $\,\gamma_{21}$}
                 & \footnotesize $pq(2rw+p-q-r)$ & \footnotesize $qr(w(p\!-\!r)\!-\!p+q+r)$ & \footnotesize $r(w(p\!-\!q)\!-\!p+q+r)$
\\
\hline
\small $\!\! $  & \!\!\! \mbox{\small $\,\alpha_{22}$}
                & \footnotesize $\!-pq\!-\!(p\!+\!q)wr$ & \footnotesize $r(w(r-q)+q)$ & \small $r$
\\
 $\Pi_{22}$ & \!\!\! \mbox{\small $\,\beta_{22}$}
                       & \footnotesize $-pq(w+1)$ & \footnotesize $q(r-pw)$ & \footnotesize $w(q\!-\!p\!-\!r)\!+\!r$
\\
\scriptsize $x\geq 0$ & \!\! \mbox{\small $\,\gamma_{22}$}
                      & \footnotesize $pq(2rw+p-q-r)$ & \footnotesize $qr(w(p\!-\!r)\!-\!p+q+r)$ & \footnotesize $r(w(p\!-\!q)\!-\!p+q+r)$
\\
\hline
\end{tabular}
\end{center}
\begin{center}
  \textit{Table 2: The \textsc{Erd\" os-Mordell} inequality in the taxicab plane geometry for case $2^{\circ}$}
\end{center}

\bigskip

Let us notice that the \textsc{Erd\" os-Mordell} inequality reduces to a problem of the positivity of the linear function
$$
f_{ij}(x,y) = \alpha_{ij} x + \beta_{ij} y + \gamma_{ij} \geq 0,
$$
for some choice of interior points $(x,y)$ of a triangle, for concretely defined values of parameters $\alpha_{ij}, \beta_{ij}$ and $\gamma_{ij}$ given by the above tables.
The problem of determining the minimum and maximum of linear functions $f_{ij}(x,y)$ reduces down to the determining of the minimum and maximum in the vertices of the considered triangles, according to \cite{Kasana_Kuma_2004}. Given that, it is enough to consider the cases of the minima and maxima of linear functions $f_{ij}(x,y)$ in vertices of
$\triangle \, ABD$ and $\triangle \, BCD$ for $A(0,r)$, $B(p,0)$, $C(q,0)$ and $D(p,\frac{r}{q}(q-p))$ when $i=1$ and in vertices of $\triangle \, ABO$ and $\triangle \, ACO$
for $A(0,r)$, $B(p,0)$, $C(q,0)$ and $O(0,0)$ when $i=2$.

\medskip
The following statements hold:

\begin{statement}
Let $A(0,r)\in \left[\Pi_{11}\right]$. If the inequality $(\ref{E8})$ holds for $A(0,r)$, then the following conclusion holds for the weight coefficient $w\!:$
  \begin{equation}
    \label{E14}
\langle \textbf{\textit{a}} \rangle  \, 0 \!<\! r \!\leq\! p\! <\! q \; \vee \;
\langle \textbf{\textit{b}} \rangle  \, 0 \!\leq\! p \!<\! r \!\leq\! q  \; \vee \;
\langle \textbf{\textit{c}} \rangle \,  0 \!\leq\! p \!<\! q \!<\! r
\;\Longrightarrow\;
      w \leq 2+ \displaystyle\frac{p+q}{r}.
  \end{equation}
\end{statement}
\begin{proof} From Table 1: \\[1.0 ex]
$\langle \textbf{\textit{a}} \rangle$ \, By substituting coordinates $x = 0$ and $y = r$ into $f_{11}(x,y)=\alpha_{11} x + \beta_{11} y + \gamma_{11}$
the following is obtained:
$$
\begin{array}{ccl}
f_{11}(0,r) \geq 0  \!&\! \;\Longleftrightarrow\;  \!&\! ((p-q)wr-pq)\cdot 0 -pq(w-1)\cdot r + pq(p+q+r) \geq 0                                             \\[2.5 ex]
                    \!&\! \;\Longleftrightarrow\;  \!&\! -pq(w-1) \cdot r+pq(p+q+r)\geq 0                                                                   \\[2.5 ex]
                    \!&\! \;\mathop{\Longrightarrow}\limits_{\mbox{\tiny $pq > 0$}}\;  \!&\! -wr + p+q+2r \geq 0                                            \\[1.5 ex]
                    \!&\! \;\mathop{\Longrightarrow}\limits_{\mbox{\tiny $r > 0$}}\;
                                                   \!&\! w \leq 2+ \displaystyle\frac{p+q}{r};
\end{array}
$$
$\langle \textbf{\textit{b}} \rangle$\, $q(r-pw) \cdot r+qr(w(p-r)+p+q+r)\ge 0$, from which follows $w \leq 2+ \displaystyle\frac{p+q}{r}$;

\medskip
\noindent
$\langle \textbf{\textit{c}} \rangle$\, $w(q\!-\!p\!-\!r)\cdot r+r(w(p\!-\!q)\!+\!p\!+\!q\!+\!r)\!+\!r\ge 0$, from which follows $w \leq 2+ \displaystyle\frac{p\!+\!q}{r}$.

\hfill $\Box$ \end{proof}

\begin{statement}
Let $A(0,r)\in \left[\Pi_{22}\right]$. If the inequality  $(\ref{E8})$ holds for $A(0,r)$, then the following conclusion holds for the weight coefficient $w\!:$
  \begin{equation}
    \label{E21}
 \langle \textbf{\textit{a}} \rangle  \, 0 \!<\! r \!\leq\! -p\! \leq \! q \; \vee \;
\langle \textbf{\textit{b}} \rangle  \, 0 \!<\! -p \!\leq\! r \!\leq\! q  \; \vee \;
\langle \textbf{\textit{c}} \rangle \,  0 \!<\! -p \!\leq\! q \!<\! r
\;\Longrightarrow\;
        w \leq 2+ \displaystyle\frac{q-p}{r}.
  \end{equation}
\end{statement}
\begin{proof} By Table 2: \\[1.0 ex]
$\langle \textbf{\textit{a}} \rangle$\, $-pq(w+1)r+pq(2rw+p-q-r)\ge 0$, from which follows $w \leq 2+ \displaystyle\frac{q-p}{r}$;
\\
$\langle \textbf{\textit{b}} \rangle$\, $q(r-pw)r+qr(w(p-r)-p+q+r)\ge 0$, from which follows $w \leq 2+ \displaystyle\frac{q-p}{r}$;
\\
$\langle \textbf{\textit{c}} \rangle$\, $(w(q\!-\!p\!-\!r)+r)r+r(w(p\!-\!q)\!-\!p\!+\!q\!+\!r) \ge 0$, from which follows $w \leq 2+ \displaystyle\frac{q-p}{r}$.

\hfill $\Box$ \end{proof}

\begin{statement}
Let $B(p,0)\in \left[\Pi_{12}\right]$. If the inequality  $(\ref{E8})$ holds for $B(p,0)$, then the following conclusions hold for the weight coefficient $w\!:$
\begin{equation}
    \label{E15}
\langle \textbf{\textit{a}} \rangle  \, 0 \!<\! r \!\leq\! p\! <\! q \; \vee \;
\langle \textbf{\textit{b}} \rangle  \, 0 \!\leq\! p \!<\! r \!\leq\! q
\;\Longrightarrow\;
      w\leq 1+ \displaystyle\frac{q^2+pr}{r(q-p)};
\end{equation}
\begin{equation}
    \label{E16}
\langle \textbf{\textit{c}} \rangle \,  0 \!\leq\! p \!<\! q \!<\! r
\;\Longrightarrow\;
     w\leq 1+ \displaystyle\frac{p+r}{q-p}.
\end{equation}
\end{statement}

\begin{proof} By Table 1: \\[1.0 ex]
$\langle \textbf{\textit{a}} \rangle$\, $((p-q)wr+pq)p+pq(-p+q+r)\ge 0$, from which follows  $w\leq 1+ \displaystyle\frac{q^2+pr}{r(q-p)}$;
\\
$\langle \textbf{\textit{b}} \rangle$\, $r(w(r\!-\!q)\!+\!q)p\!+\!qr(w(p\!-\!r)\!-\!p\!+\!q\!+\!r) \!\ge\! 0$, from which follows
$w \!\leq\! 1 \!+\! \displaystyle\frac{q^2\!+\!pr}{r(q\!-\!p)}$;
\\
$\langle \textbf{\textit{c}} \rangle$\, $rp+r(w(p-q)-p+q+r)\ge 0$, from which follows  $w\leq 1+ \displaystyle\frac{p+r}{q-p}$.
\hfill $\Box$ \end{proof}

\begin{statement}
Let $B(p,0)\in \left[\Pi_{21}\right]$. If the inequality  $(\ref{E8})$ holds for $B(p,0)$, then the following conclusions hold for the weight coefficient $w\!:$
\begin{equation}
    \label{E17}
\langle \textbf{\textit{a}} \rangle  \, 0 \!<\! r \!\leq\! -p\! \leq\! q \; \vee \;
\langle \textbf{\textit{b}} \rangle  \, 0 \!<\! -p \!\leq\! r \!\leq\! q
\;\Longrightarrow\;
       w \leq \displaystyle\frac{q}{r}\left(1+ \frac{r-p}{q-p}\right);
\end{equation}
\begin{equation}
    \label{E18}
\langle \textbf{\textit{c}} \rangle \,  0 \!<\! -p \!\leq\! q \!<\! r
\;\Longrightarrow\;
     w\leq 1+ \displaystyle\frac{r-p}{q-p}.
\end{equation}
\end{statement}
\begin{proof} By Table 2: \\[1.0 ex]
$\langle \textbf{\textit{a}} \rangle$\, $(pq-(p+q)wr)p+pq(2rw+p-q-r)\ge 0$, from which follows $w \leq \displaystyle\frac{q}{r}\left(\!1+ \frac{r-p}{q-p}\!\right)$;
\\
$\langle \textbf{\textit{b}} \rangle$\, $r(w(r-q)-q)p+qr(w(p-r)-p+q+r)\ge 0$, from which follows $w \leq \displaystyle\frac{q}{r}\left(\!1+ \frac{r-p}{q-p}\!\right)$;
\\
$\langle \textbf{\textit{c}} \rangle$\, $-rp+r(w(p-q)-p+q+r)\ge 0$, from which follows $w\leq 1+ \displaystyle\frac{r-p}{q-p}$.
\hfill $\Box$ \end{proof}

\begin{statement}
Let $C(q,0)\in \left[\Pi_{12}\right]$. If the inequality $(\ref{E8})$ holds for $C(q,0)$, then the following conclusions hold for the weight coefficient $w\!:$
 \begin{equation}
    \label{E19}
\langle \textbf{\textit{a}} \rangle  \, 0 \!<\! r \!\leq\! p\! <\! q
\;\Longrightarrow\;
      w\leq \frac{p}{r}\left(1+ \displaystyle\frac{q+r}{q-p}\right);
  \end{equation}
  \begin{equation}
    \label{E20}
\langle \textbf{\textit{b}} \rangle  \, 0 \!\leq\! p \!<\! r \!\leq\! q  \; \vee \;
\langle \textbf{\textit{c}} \rangle \,  0 \!\leq\! p \!<\! q \!<\! r
\;\Longrightarrow\;
    w\leq 1+ \displaystyle\frac{q+r}{q-p}.
  \end{equation}
\end{statement}

\begin{proof} By Table 1: \\
$\langle \textbf{\textit{a}} \rangle$\, $((p-q)wr+pq)q+pq(-p+q+r)\ge 0$,
from which follows $w\leq \displaystyle\frac{p}{r}\left(1+ \displaystyle\frac{q+r}{q-p}\right)$;
\\
$\langle \textbf{\textit{b}} \rangle$\, $r(w(r\!-\!q)\!+\!q)q\!+\!qr(w(p\!-\!r)\!-\!p\!+\!q\!+\!r) \!\ge\! 0$, from which follows
$w \!\leq\! 1 \!+\! \displaystyle\frac{q\!+\!r}{q\!-\!p}$;
\\
$\langle \textbf{\textit{c}} \rangle$\, $rq+r(w(p-q)-p+q+r)\ge 0$, from which follows $w\leq 1+ \displaystyle\frac{q+r}{q-p}$.
\hfill $\Box$ \end{proof}

\begin{statement}
Let $C(q,0)\in \left[\Pi_{22}\right]$. If the inequality $(\ref{E8})$ holds for $C(q,0)$, then the following conclusions hold for the weight coefficient $w\!:$
  \begin{equation}
    \label{E21}
\langle \textbf{\textit{a}} \rangle  \, 0 \!<\! r \!\leq\! -p\! \leq\! q
\;\Longrightarrow\;
      w \leq \displaystyle\frac{-p}{r}\left(1+ \displaystyle\frac{q+r}{q-p}\right);
  \end{equation}
  \begin{equation}
    \label{E22}
\langle \textbf{\textit{b}} \rangle  \, 0 \!<\! -p \!\leq\! r \!\leq\! q  \; \vee \;
\langle \textbf{\textit{c}} \rangle \,  0 \!<\! -p \!\leq\! q \!<\! r
\;\Longrightarrow\;
    w\leq 1+ \displaystyle\frac{q+r}{q-p}.
  \end{equation}
\end{statement}
\begin{proof} By Table 2: \\[1.0 ex]
$\langle \textbf{\textit{a}} \rangle$\, $(-pq\!-\!(p\!+\!q)wr)q\!+\!pq(2rw\!+\!p\!-\!q\!-\!r) \!\ge\! 0$,
from which follows $w \!\leq\! \displaystyle\frac{-p}{r}\left(\!1\!+\!\displaystyle\frac{q\!+\!r}{q\!-\!p}\!\right)$;
\\
$\langle \textbf{\textit{b}} \rangle$\, $r(w(r\!-\!q)\!+\!q)q\!+\!qr(w(p\!-\!r)\!-\!p\!+\!q\!+\!r) \!\ge\! 0$, from which follows
$w \!\leq\! 1 \!+\! \displaystyle\frac{q\!+\!r}{q\!-\!p}$;
\\
$\langle \textbf{\textit{c}} \rangle$\, $rq+r(w(p-q)-p+q+r)\ge 0$, from which follows $w\leq 1+ \displaystyle\frac{q+r}{q-p}$.
\hfill $\Box$ \end{proof}

\begin{statement}
Let $D(p,\frac{r}{q}(q-p))\in \left[\Pi_{12}\right]$. If the inequality $(\ref{E8})$ holds for $D(p,\frac{r}{q}(q-p))$,
then the following conclusions hold for the weight coefficient $w\!:$
 \begin{equation}
    \label{E23}
\langle \textbf{\textit{a}} \rangle  \, 0 \!<\! r \!\leq\! p\! <\! q
\;\Longrightarrow\;
      w\leq 1+ \displaystyle\frac{q^2+pr}{2r(q-p)};
  \end{equation}
  \begin{equation}
    \label{E24}
\langle \textbf{\textit{b}} \rangle  \, 0 \!\leq\! p \!<\! r \!\leq\! q  \; \vee \;
\langle \textbf{\textit{c}} \rangle \,  0 \!\leq\! p \!<\! q \!<\! r
\;\Longrightarrow\;
    w\leq 1+ \displaystyle\frac{q-p}{r+p}+ \displaystyle\frac{q}{q-p}.
  \end{equation}
\end{statement}
\begin{proof} By Table 1: \\[1.0 ex]
$\langle \textbf{\textit{a}} \rangle$\, $((p-q)wr+pq)p-pq(w-1)\displaystyle\frac{r}{q}(q-p)+pq(-p+q+r)\ge 0$, \\
from which follows $w\leq 1+ \displaystyle\frac{q^2+pr}{2r(q-p)}$;
\\
$\langle \textbf{\textit{b}} \rangle$\, $r(w(r-q)+q)p+q(r-pw)\displaystyle\frac{r}{q}(q-p)+qr(w(p-r)-p+q+r)\ge 0$, \\
from which follows $w\leq 1+ \displaystyle\frac{q-p}{r+p}+ \displaystyle\frac{q}{q-p}$;
\\
$\langle \textbf{\textit{c}} \rangle$\, $rp+(w(q-p-r)+r)\displaystyle\frac{r}{q}(q-p)+r(w(p-q)-p+q+r)\ge 0$, \\
from which follows $w\leq 1+ \displaystyle\frac{q-p}{r+p}+ \displaystyle\frac{q}{q-p}$.
\hfill $\Box$ \end{proof}

\begin{statement}
Let $O(0,0)\in \left[\Pi_{22}\right]$. If the inequality $(\ref{E8})$ holds for $O(0,0)$, then the following conclusions hold for the weight coefficient $w\!:$
 \begin{equation}
    \label{E25}
\langle \textbf{\textit{a}} \rangle  \, 0 \!<\! r \!\leq\! -p\! \leq\! q
\;\Longrightarrow\;
            w\leq  \displaystyle\frac{1}{2}+ \displaystyle\frac{q-p}{2r};
  \end{equation}
 \begin{equation}
    \label{E26}
\langle \textbf{\textit{b}} \rangle  \, 0 \!<\! -p \!\leq\! r \!\leq\! q
\;\Longrightarrow\;
    w\leq  1+ \displaystyle\frac{q}{r-p};
  \end{equation}
  \begin{equation}
    \label{E27}
\langle \textbf{\textit{c}} \rangle \,  0 \!<\! -p \!\leq\! q \!<\! r
\;\Longrightarrow\;
    w\leq 1+ \displaystyle\frac{r}{q-p}.
  \end{equation}
\end{statement}
\begin{proof} By Table 2: \\[1.0 ex]
$\langle \textbf{\textit{a}} \rangle$\,$pq(2rw+p-q-r)\ge 0$, from which follows $w\leq \displaystyle\frac{1}{2}+ \displaystyle\frac{q-p}{2r}$;
\\
$\langle \textbf{\textit{b}} \rangle$\, $qr(w(p-r)-p+q+r)\ge 0$, from which follows $w\leq  1+ \displaystyle\frac{q}{r-p}$;
\\
$\langle \textbf{\textit{c}} \rangle$\, $r(w(p-q)-p+q+r)\ge 0$, from which follows $w\leq 1+ \displaystyle\frac{r}{q-p}$.
\hfill $\Box$ \end{proof}

\bigskip
Let the positions of points $B$ and $C$ be given. Then, let us consider the positions of point $A(0, r)$ in the concrete cases
$\langle \textbf{\textit{a}} \rangle$,
$\langle \textbf{\textit{b}} \rangle$,
$\langle \textbf{\textit{c}} \rangle$
which were considered in Statements 1--8. Through the aforementioned Statements the functions of upper bounds $\omega$ for the weight coefficient $w$ were obtained:
$$
w \leq \omega(p,q,r).
$$
Our goal is to, for the functions $\omega(p,q,r)$, dependent on concrete subcases $\langle \mbox{\boldmath $\theta$} \rangle$, where $\mbox{\boldmath $\theta$} \in \{\textbf{\textit{a}}, \textbf{\textit{b}}, \textbf{\textit{c}}\}$, find the values:
\begin{equation}
\mathfrak{M} = \inf\{\omega(p,q,r) \, | \, \langle \mbox{\boldmath $\theta$} \rangle  \}.
\end{equation}
In this way, the \textsc{Erd\" os-Mordell} inequality (\ref{E7}) holds for $w=\mathfrak{M}$ for all interior points of $\triangle \, ABC$. If $\mathfrak{M}$ is a minimum in this area, then an equality is also possible in (\ref{E7}).

\bigskip
\noindent
{\bf 2.1$\;\;$Determining value of {\boldmath $\mathfrak{M}$} by areas}

\bigskip
In this section of the work, the values of $\mathfrak{M}$ by areas of $\triangle ABC$ are determined in dependence on cases $\langle \mbox{\boldmath $\theta$} \rangle$,
where $\mbox{\boldmath $\theta$} \in \{\textbf{\textit{a}}, \textbf{\textit{b}}, \textbf{\textit{c}}\}$.


\bigskip
The following three propositions are obtained on the basis of Statement 1.

\begin{proposition}

Let $A(0,r)\in \left[\Pi_{11}\right]$. If the inequality $(\ref{E8})$ holds for  $A(0,r)$, then the following conclusion holds for the weight coefficient $w\!:$
  \begin{equation}
\langle \textbf{\textit{a}} \rangle  \; 0 \!<\! r \!\leq\! p\! <\! q
\;\Longrightarrow\;
      w \leq \omega(p,q,r) = 2+ \displaystyle\frac{p+q}{r}
  \end{equation}
and in that case
  \begin{equation}
  \omega(p,q,r) \in (\mathfrak{M}, \infty) \;\;\mbox{and}\;\; \mathfrak{M}=4.
  \end{equation}
 \end{proposition} \begin{proof}
Let us consider $\langle \textbf{\textit{a}} \rangle  \; 0 \!<\! r \!\leq\! p\! <\! q $. Then, we notice the following expression holds:
$$
\omega(p,q,r)
=
2 + \displaystyle\frac{p+q}{r}
\geq
2 + \displaystyle\frac{p+q}{p}
=
3 + \displaystyle\frac{q}{p}
>
4
 \;\Longrightarrow\; \mathfrak{M} = 4.
$$
The above conclusion is correct because the real number $\displaystyle\frac{q}{p}$ fulfills $\displaystyle\frac{q}{p} > 1$ and it is possible to choose a number $\displaystyle\frac{q}{p}$ such that it is arbitrarily close to $1$.
\hfill $\Box$ \end{proof}

\begin{proposition}
Let $A(0,r)\in \left[\Pi_{11}\right]$. If the inequality $(\ref{E8})$ holds for  $A(0,r)$, then the following conclusion holds for the weight coefficient $w\!:$
  \begin{equation}
\langle \textbf{\textit{b}} \rangle  \; 0  \!\leq\! p\! < \! r \! \leq \! q
\;\Longrightarrow\;
      w \leq \omega(p,q,r) = 2+ \displaystyle\frac{p+q}{r}
  \end{equation}
and in that case
  \begin{equation}
  \omega(p,q,r) \in [\mathfrak{M}, \infty) \;\;\mbox{and}\;\; \mathfrak{M}=3.
  \end{equation}
 \end{proposition} \begin{proof}
Let us consider $\langle \textbf{\textit{b}} \rangle  \; 0  \!\leq\! p\! < \! r \! \leq \! q$. Then, we notice the following expression holds:
$$
\omega(p,q,r)
=
2 + \displaystyle\frac{p+q}{r}
\geq
2 + \displaystyle\frac{p+q}{q}
=
3 + \displaystyle\frac{p}{q}
\geq
3
 \;\Longrightarrow\; \mathfrak{M} = 3.
$$
\hfill $\Box$ \end{proof}

\begin{proposition}
Let $A(0,r)\in \left[\Pi_{11}\right]$. If the inequality $(\ref{E8})$ holds for $A(0,r)$, then the following conclusion holds for the weight coefficient $w\!:$
  \begin{equation}
\langle \textbf{\textit{c}} \rangle  \; 0  \!\leq\! p\! < \! q \! < \! r
\;\Longrightarrow\;
      w \leq \omega(p,q,r) = 2+ \displaystyle\frac{p+q}{r}
  \end{equation}
and in that case
  \begin{equation}
  \omega(p,q,r) \in (\mathfrak{M}, \infty) \;\;\mbox{and}\;\; \mathfrak{M}=2.
  \end{equation}
 \end{proposition} \begin{proof}
Let us consider $\langle \textbf{\textit{c}} \rangle  \; 0  \!\leq\! p\! < \! q \! < \! r$. Then, we notice the following expression holds:
$$
\omega(p,q,r)
=
2 + \displaystyle\frac{p+q}{r}
>
2
 \;\Longrightarrow\; \mathfrak{M} = 2.
$$
\hfill $\Box$ \end{proof}


The following three propositions are obtained on the basis of Statement 2.

\begin{proposition}
Let $A(0,r)\in \left[\Pi_{22}\right]$. If the inequality $(\ref{E8})$ holds for $A(0,r)$, then the following conclusion holds for the weight coefficient $w\!:$
  \begin{equation}
\langle \textbf{\textit{a}} \rangle  \; 0 \!<\! r \!\leq\! -p\! \leq \! q
\;\Longrightarrow\;
      w \leq \omega(p,q,r) = 2+ \displaystyle\frac{q-p}{r}
  \end{equation}
and in that case
  \begin{equation}
  \omega(p,q,r) \in [\mathfrak{M}, \infty) \;\;\mbox{and}\;\; \mathfrak{M}=4.
  \end{equation}
 \end{proposition}
\begin{proof}
Let us consider $\langle \textbf{\textit{a}} \rangle  \; 0 \!<\! r \!\leq\! -p\! \leq \! q $. Then, we notice the following expression holds:
$$
\omega(p,q,r)
=
2 + \displaystyle\frac{q-p}{r}
\geq
2 + \displaystyle\frac{q-p}{-p}
=
3 + \displaystyle\frac{q}{-p}
\geq
4
 \;\Longrightarrow\; \mathfrak{M} = 4.
$$
\hfill $\Box$ \end{proof}

\begin{proposition}
Let $A(0,r)\in \left[\Pi_{22}\right]$. If the inequality $(\ref{E8})$ holds for $A(0,r)$, then the following conclusion holds for the weight coefficient $w\!:$
  \begin{equation}
\langle \textbf{\textit{b}} \rangle  \; 0 \! < \! -p \! \leq \! r\! \leq \! q
\;\Longrightarrow\;
      w \leq \omega(p,q,r) = 2 + \displaystyle\frac{q-p}{r}
  \end{equation}
and in that case
  \begin{equation}
  \omega(p,q,r) \in (\mathfrak{M}, \infty) \;\;\mbox{and}\;\; \mathfrak{M}=3.
  \end{equation}
 \end{proposition}
\begin{proof}
Let us consider $\langle \textbf{\textit{b}} \rangle  \; 0 \! < \! -p \! \leq \! r \! \leq \! q $. Then, we notice the following expression holds:
$$
\omega(p,q,r)
=
2 + \displaystyle\frac{q-p}{r}
\geq
2 + \displaystyle\frac{q-p}{q}
=
3 + \displaystyle\frac{-p}{q}
>
3
 \;\Longrightarrow\; \mathfrak{M} = 3.
$$
\hfill $\Box$ \end{proof}

\begin{proposition}
Let $A(0,r)\in \left[\Pi_{22}\right]$. If the inequality $(\ref{E8})$ holds for $A(0,r)$, then the following conclusion holds for the weight coefficient $w\!:$
  \begin{equation}
\langle \textbf{\textit{c}} \rangle  \; 0 \! < \! -p \! \leq \! q\! < \! r
\;\Longrightarrow\;
      w \leq \omega(p,q,r) = 2 + \displaystyle\frac{q-p}{r}
  \end{equation}
and in that case
  \begin{equation}
  \omega(p,q,r) \in (\mathfrak{M}, \infty) \;\;\mbox{and}\;\; \mathfrak{M}=2.
  \end{equation}
 \end{proposition}
 \begin{proof}
Let us consider $\langle \textbf{\textit{c}} \rangle  \; 0 \! < \! -p \! \leq \! q \! < \! r $. Then, we notice the following expression holds:
$$
\omega(p,q,r)
=
2 + \displaystyle\frac{q-p}{r}
>
2
\;\Longrightarrow\; \mathfrak{M} = 2.
$$
\hfill $\Box$ \end{proof}


Similar to previous propositions, the following three propositions are obtained from Statement 3.

\begin{proposition}
Let $B(p,0)\in \left[\Pi_{12}\right]$. If the inequality $(\ref{E8})$ holds for $B(p,0)$, then the following conclusion holds for the weight coefficient $w\!:$
  \begin{equation}
\langle \textbf{\textit{a}} \rangle  \; 0 \!<\! r \!\leq\! p\! <\! q
\;\Longrightarrow\;
      w \leq \omega(p,q,r) = 1 + \displaystyle\frac{q^2}{r(q-p)} + \displaystyle\frac{p}{q-p}
  \end{equation}
and in that case
  \begin{equation}
  \omega(p,q,r) \in [\mathfrak{M}, \infty) \;\;\mbox{and}\;\; \mathfrak{M}=3+2\sqrt{2}.
  \end{equation}
 \end{proposition}
\begin{proof}
Let us consider $\langle \textbf{\textit{a}} \rangle  \; 0 \!<\! r \!\leq\! p\! <\! q $. Then, we notice the following expression holds:
$$
\begin{array}{rclc}
\omega(p,q,r)
\!&\!=\!&\!
1 + \displaystyle\frac{q^2}{r(q-p)} + \displaystyle\frac{p}{q-p}                &                                      \\[2.5 ex]
\!&\!\geq\!&\!
1 + \displaystyle\frac{q^2}{p(q-p)} + \displaystyle\frac{p}{q-p}                &                                      \\[2.5 ex]
\!&\! = \!&\!
3 + \displaystyle\frac{2p}{q-p} + \displaystyle\frac{q-p}{p} \geq 3 + 2\sqrt{2} & \;\Longrightarrow\; \mathfrak{M} = 3+2\sqrt{2},
\end{array}
$$

\noindent
because $t = \displaystyle\frac{p}{q - p} > 0$ holds $2 \,t + \mbox{\small $\displaystyle\frac{1}{t}$} \geq 2 \sqrt{2}$.
\hfill $\Box$ \end{proof}

\begin{proposition}
Let $B(p,0)\in \left[\Pi_{12}\right]$.   If inequality $(\ref{E8})$ holds for $B(p,0)$, then the following conclusion holds for the weight coefficient $w\!:$
  \begin{equation}
\langle \textbf{\textit{b}} \rangle  \; 0 \! \leq \! p \! < \! r\! \leq \! q
\;\Longrightarrow\;
      w \leq \omega(p,q,r) = 1 + \displaystyle\frac{q^2}{r(q-p)} + \displaystyle\frac{p}{q-p}
  \end{equation}
and in that case
  \begin{equation}
  \omega(p,q,r) \in [\mathfrak{M}, \infty) \;\;\mbox{and}\;\; \mathfrak{M}=2.
  \end{equation}
 \end{proposition}
\begin{proof}
Let us consider $\langle \textbf{\textit{b}} \rangle  \; 0 \! \leq \! p \! < \! r\! \leq \! q$. Then, we notice the following expression holds:
$$
\begin{array}{rclc}
\omega(p,q,r)
\!&\!=\!&\!
1 + \displaystyle\frac{q^2}{r(q-p)} + \displaystyle\frac{p}{q-p}   &                                      \\[2.5 ex]
\!&\!\geq\!&\!
1 + \displaystyle\frac{q^2}{q(q-p)} + \displaystyle\frac{p}{q-p}   &                                      \\[2.5 ex]
\!&\! = \!&\!
2 + \displaystyle\frac{2p}{q-p} \geq 2                             & \;\Longrightarrow\; \mathfrak{M} = 2.
\end{array}
$$
\hfill $\Box$ \end{proof}

\begin{proposition}
Let $B(p,0)\in \left[\Pi_{12}\right]$.  If the inequality $(\ref{E8})$ holds for $B(p,0)$, then the following conclusion holds for the weight coefficient $w\!:$
  \begin{equation}
\langle \textbf{\textit{c}} \rangle  \; 0 \! \leq \! p \! < \! q\! < \! r
\;\Longrightarrow\;
      w \leq \omega(p,q,r) = 1 + \displaystyle\frac{p+r}{q-p}
  \end{equation}
and in that case
  \begin{equation}
  \omega(p,q,r) \in (\mathfrak{M}, \infty) \;\;\mbox{and}\;\; \mathfrak{M}=2.
  \end{equation}
 \end{proposition}
\begin{proof}
Let us consider $\langle \textbf{\textit{c}} \rangle  \; 0 \! \leq \! p \! < \! q\! < \! r$. Then, we notice the following expression holds:
$$
\begin{array}{rclc}
\omega(p,q,r)
\!&\!=\!&\!
1 + \displaystyle\frac{p+r}{q-p}       &                                      \\[2.5 ex]
\!&\!>\!&\!
1 + \displaystyle\frac{p+q}{q-p}       &                                      \\[2.5 ex]
\!&\! = \!&\!
2 + \displaystyle\frac{2p}{q-p} \geq 2 & \;\Longrightarrow\; \mathfrak{M} = 2.
\end{array}
$$
\hfill $\Box$ \end{proof}


The following three propositions are obtained on the basis of Statement 4.

\begin{proposition}
Let $B(p,0)\in \left[\Pi_{21}\right]$. If the inequality $(\ref{E8})$ holds for $B(p,0)$, then the following conclusion holds for the weight coefficient $w\!:$
  \begin{equation}
\langle \textbf{\textit{a}} \rangle  \; 0 \!<\! r \!\leq\! -p\! \leq \! q
\;\Longrightarrow\;
      w \leq \omega(p,q,r) = \displaystyle\frac{q}{r}\left(1 + \displaystyle\frac{r-p}{q-p}\right)
  \end{equation}
and in that case
  \begin{equation}
  \omega(p,q,r) \in [\mathfrak{M}, \infty) \;\;\mbox{and}\;\; \mathfrak{M}=2.
  \end{equation}
 \end{proposition}
\begin{proof}
Let us consider $\langle \textbf{\textit{a}} \rangle  \; 0 \!<\! r \!\leq\! -p\! \leq\! q $. Then, we notice the following expression holds:
$$
\begin{array}{rclc}
\omega(p,q,r)
\!&\! = \!&\!
\displaystyle\frac{q}{r}\left(1 + \displaystyle\frac{r-p}{q-p}\right)           &                                      \\[2.5 ex]
\!&\! \geq \!&\!
\displaystyle\frac{q}{r}\left(1 + \displaystyle\frac{r-p}{2q}\right)            &                                      \\[2.5 ex]
\!&\! \geq \!&\!
\displaystyle\frac{q}{r}\left(1 + \displaystyle\frac{2r}{2q}\right) = \displaystyle\frac{q}{r} + 1 \geq 2 & \;\Longrightarrow\; \mathfrak{M} = 2.
\end{array}
$$
\hfill $\Box$ \end{proof}

\begin{proposition}
Let $B(p,0)\in \left[\Pi_{21}\right]$. If the inequality $(\ref{E8})$ holds for $B(p,0)$, then the following conclusion holds for the weight coefficient $w\!:$
  \begin{equation}
\langle \textbf{\textit{b}} \rangle  \; 0 \! < \! -p \leq r \leq q
\;\Longrightarrow\;
      w \leq \omega(p,q,r) = \displaystyle\frac{q}{r}\left(1 + \displaystyle\frac{r-p}{q-p}\right)
  \end{equation}
and in that case
  \begin{equation}
  \omega(p,q,r) \in [\mathfrak{M}, \infty) \;\;\mbox{and}\;\; \mathfrak{M}=2.
  \end{equation}
 \end{proposition}
\begin{proof}
Let us consider $\langle \textbf{\textit{b}} \rangle  \; 0 \!<\! -p \leq r \leq q$. Then, we notice the following expression holds:
$$
\omega(p,q,r) =
\displaystyle\frac{q}{r} + \displaystyle\frac{q}{r} \, \displaystyle\frac{r-p}{q-p}
\geq
1 + \displaystyle\frac{q}{r} \, \displaystyle\frac{r-p}{q-p} \geq 2 \;\Longrightarrow\; \mathfrak{M} = 2,
$$
because $q(r+(-p)) \geq r (q+(-p)) \;\Longleftrightarrow\; q \geq r$.
\hfill $\Box$ \end{proof}

\begin{proposition}
Let $B(p,0)\in \left[\Pi_{21}\right]$. If the inequality $(\ref{E8})$ holds for $B(p,0)$, then the following conclusion holds for the weight coefficient $w\!:$
  \begin{equation}
\langle \textbf{\textit{c}} \rangle  \; 0 \!<\! -p \leq q < r
\;\Longrightarrow\;
      w \leq \omega(p,q,r) = 1 + \displaystyle\frac{r-p}{q-p}
  \end{equation}
and in that case
  \begin{equation}
  \omega(p,q,r) \in (\mathfrak{M}, \infty) \;\;\mbox{and}\;\; \mathfrak{M}=2.
  \end{equation}
 \end{proposition}
\begin{proof}
Let us consider $\langle \textbf{\textit{c}} \rangle  \; 0 \!<\! -p \leq q < r$. Then, we notice the following expression holds:
$$
\omega(p,q,r)
=
1 + \displaystyle\frac{r-p}{q-p} > 1 + \displaystyle\frac{q-p}{q-p} = 2  \;\;\Longrightarrow\;\;  \mathfrak{M} = 2.
$$
\hfill $\Box$ \end{proof}


The following three propositions are obtained on the basis of Statement 5.

\begin{proposition}
Let $C(q,0)\in \left[\Pi_{12}\right]$. If the inequality $(\ref{E8})$ holds for $C(q,0)$, then the following conclusion holds for the weight coefficient $w\!:$
  \begin{equation}
\langle \textbf{\textit{a}} \rangle  \; 0 \!<\! r \!\leq\! p\! <\! q
\;\Longrightarrow\;
      w \leq \omega(p,q,r) = \displaystyle\frac{p}{r} \left(1 + \displaystyle\frac{q+r}{q-p}\right)
  \end{equation}
and in that case
  \begin{equation}
  \omega(p,q,r) \in (\mathfrak{M}, \infty) \;\;\mbox{and}\;\; \mathfrak{M}=2.
  \end{equation}
 \end{proposition}
\begin{proof}
Let us consider $\langle \textbf{\textit{a}} \rangle  \; 0 \!<\! r \!\leq\! p\! <\! q $. Then, we notice the following expression holds:
$$
\begin{array}{rclc}
\omega(p,q,r)
\!&\!=\!&\!
\displaystyle\frac{p}{r} \left(1 + \displaystyle\frac{q+r}{q-p}\right)                  &                                      \\[2.5 ex]
\!&\!\geq\!&\!
1 + \displaystyle\frac{q+r}{q-p}                                                        &                                      \\[2.5 ex]
\!&\!=\!&\!
\displaystyle\frac{2q - 2p + r + p}{q-p} \,>\,2 + \displaystyle\frac{r+p}{q-p} \,>\, 2  & \;\Longrightarrow\;\; \mathfrak{M} = 2.
\end{array}
$$
\hfill $\Box$ \end{proof}

\begin{proposition}
Let $C(q,0)\in \left[\Pi_{12}\right]$. If the inequality $(\ref{E8})$ holds for $C(q,0)$, then the following conclusion holds for the weight coefficient $w\!:$
  \begin{equation}
\langle \textbf{\textit{b}} \rangle  \; 0 \!\leq\! p \!<\! r \! \leq \! q
\;\Longrightarrow\;
      w \leq \omega(p,q,r) = 1 + \displaystyle\frac{q+r}{q-p}
  \end{equation}
and in that case
  \begin{equation}
  \omega(p,q,r) \in (\mathfrak{M}, \infty) \;\;\mbox{and}\;\; \mathfrak{M}=2.
  \end{equation}
 \end{proposition}
\begin{proof}
Let us consider  $\langle \textbf{\textit{b}} \rangle  \; 0 \!\leq\! p \!<\! r \! \leq \! q$. Then, we notice the following expression holds:
$$
\omega(p,q,r)
=
1 + \displaystyle\frac{q+r}{q-p}
>
1 + \displaystyle\frac{q+p}{q-p} \geq 2  \;\;\Longrightarrow\;\; \mathfrak{M} = 2.
$$
\hfill $\Box$ \end{proof}

\begin{proposition}
Let $C(q,0)\in \left[\Pi_{12}\right]$. If the inequality $(\ref{E8})$ holds for $C(q,0)$, then the following conclusion holds for the weight coefficient $w\!:$
  \begin{equation}
\langle \textbf{\textit{c}} \rangle  \; 0 \!\leq\! p \!<\! q \!<\! r
\;\Longrightarrow\;
      w \leq \omega(p,q,r) = 1 + \displaystyle\frac{q+r}{q-p}
  \end{equation}
and in that case
  \begin{equation}
  \omega(p,q,r) \in (\mathfrak{M}, \infty) \;\;\mbox{and}\;\; \mathfrak{M}=3.
  \end{equation}
 \end{proposition}
\begin{proof}
Let us consider  $\langle \textbf{\textit{c}} \rangle  \; 0 \!\leq\! p \!<\! q \!<\! r$. Then, we notice the following expression holds:
$$
\omega(p,q,r)
=
1 + \displaystyle\frac{q+r}{q-p}
=
2 + \displaystyle\frac{r+p}{q-p}
>
2 + \displaystyle\frac{q+p}{q-p} \,\geq\, 3 \;\;\Longrightarrow\;\; \mathfrak{M} = 3.
$$
\hfill $\Box$ \end{proof}


Similar to previous propositions, the following three propositions are obtained from Statement 6.

\begin{proposition}
Let $C(q,0)\in \left[\Pi_{22}\right]$. If the inequality $(\ref{E8})$ holds for $C(q,0)$, then the following conclusion holds for the weight coefficient $w\!:$
  \begin{equation}
\langle \textbf{\textit{a}} \rangle  \; 0 \!<\! r \!\leq\! -p\! \leq\! q
\;\Longrightarrow\;
      w \leq \omega(p,q,r) = -\displaystyle\frac{p}{r} \left(1 + \displaystyle\frac{q+r}{q-p}\right)
  \end{equation}
and in that case
  \begin{equation}
  \omega(p,q,r) \in [\mathfrak{M}, \infty) \;\;\mbox{and}\;\; \mathfrak{M}=2.
  \end{equation}
 \end{proposition}
\begin{proof}
Let us consider  $\langle \textbf{\textit{a}} \rangle  \; 0 \!<\! r \!\leq\! -p\!\leq\! q $. Then, we notice the following expression holds:
$$
\begin{array}{rclc}
\omega(p,q,r)
\!&\!=\!&\!
\displaystyle\frac{-p}{r} + \displaystyle\frac{-p}{r} \displaystyle\frac{q+r}{q-p}       &                                      \\[2.5 ex]
\!&\!\geq\!&\!
1 + \displaystyle\frac{-p}{r}\displaystyle\frac{q+r}{q-p}                                &                                      \\[2.5 ex]
\!&\!=\!&\!
1 + \displaystyle\frac{-pq + (-p)r}{rq + (-p)r} \,\geq\, 2  & \;\Longrightarrow\;\; \mathfrak{M} = 2,
\end{array}
$$
because $-pq + (-p)r \geq rq + (-p)r \; \Longleftrightarrow \; -p \geq r$.
\hfill $\Box$ \end{proof}

\begin{proposition}
Let $C(q,0)\in \left[\Pi_{22}\right]$. If the inequality $(\ref{E8})$ holds for $C(q,0)$, then the following conclusion holds for the weight coefficient $w\!:$
  \begin{equation}
\langle \textbf{\textit{b}} \rangle  \; 0 \! < \! -p \!\leq\! r \! \leq \! q
\;\Longrightarrow\;
      w \leq \omega(p,q,r) = 1 + \displaystyle\frac{q+r}{q-p}
  \end{equation}
and in that case
  \begin{equation}
  \omega(p,q,r) \in [\mathfrak{M}, \infty) \;\;\mbox{and}\;\; \mathfrak{M}=2.
  \end{equation}
 \end{proposition}
\begin{proof}
Let us consider   $\langle \textbf{\textit{b}} \rangle  \; 0 \! < \! -p \! \leq \! r \! \leq \! q$. Then, we notice the following expression holds:
$$
\omega(p,q,r)
=
1 + \displaystyle\frac{q+r}{q-p}
\geq
1 + \displaystyle\frac{q-p}{q-p} = 2  \;\;\Longrightarrow\;\; \mathfrak{M} = 2.
$$
\hfill $\Box$ \end{proof}

\begin{proposition}
Let $C(q,0)\in \left[\Pi_{22}\right]$. If the inequality $(\ref{E8})$ holds for $C(q,0)$, then the following conclusion holds for the weight coefficient $w\!:$
  \begin{equation}
\langle \textbf{\textit{c}} \rangle  \; 0 \!<\! -p \!\leq\! q \!<\! r
\;\Longrightarrow\;
      w \leq \omega(p,q,r) = 1 + \displaystyle\frac{q+r}{q-p}
  \end{equation}
and in that case
  \begin{equation}
  \omega(p,q,r) \in (\mathfrak{M}, \infty) \;\;\mbox{and}\;\; \mathfrak{M}=2.
  \end{equation}
 \end{proposition}
\begin{proof}
Let us consider $\langle \textbf{\textit{c}} \rangle  \; 0 \!<\! -p \!\leq\! q \!<\! r$. Then, we notice the following expression holds:
$$
\omega(p,q,r)
=
1 + \displaystyle\frac{q+r}{q-p}
>
1 + \displaystyle\frac{2q}{q-p}
\,\geq\, 2 \;\;\Longrightarrow\;\; \mathfrak{M} = 2,
$$
because $2q \geq q-p \;\Longleftrightarrow\; q \geq -p$.
\hfill $\Box$ \end{proof}

\bigskip


The following three propositions are obtained on the basis of Statement 7.

\begin{proposition}
Let $D{\big (}p,\frac{r}{q}(q-p){\big )} \in \left[\Pi_{12}\right]$. If the inequality $(\ref{E8})$ holds for $D{\big (}p,\frac{r}{q}(q-p){\big )} \in \left[\Pi_{12}\right]$, then the following conclusion holds for the weight coefficient $w\!:$
  \begin{equation}
\langle \textbf{\textit{a}} \rangle  \; 0 \!<\! r \!\leq\! p\! <\! q
\;\Longrightarrow\;
      w \leq \omega(p,q,r) = 1 + \displaystyle\frac{q^2+pr}{2r(q-p)}
  \end{equation}
and in that case
  \begin{equation}
  \omega(p,q,r) \in [\mathfrak{M}, \infty) \;\;\mbox{and}\;\; \mathfrak{M}=2+\sqrt{2}.
  \end{equation}
 \end{proposition}
\begin{proof}
Let us consider $\langle \textbf{\textit{a}} \rangle  \; 0 \!<\! r \!\leq\! p\! <\! q $. Then, we notice the following expression holds:
$$
\begin{array}{rclc}
\omega(p,q,r)
\!&\!=\!&\!
1 + \displaystyle\frac{q^2}{2r(q-p)} + \displaystyle\frac{p}{2(q-p)}            &                                      \\[2.5 ex]
\!&\!\geq\!&\!
1 + \displaystyle\frac{q^2}{2p(q-p)} + \displaystyle\frac{p}{2(q-p)}            &                                      \\[2.5 ex]
\!&\! = \!&\!
2 + \displaystyle\frac{q-p}{2p} + \displaystyle\frac{p}{q-p} \geq 2 + \sqrt{2} & \;\Longrightarrow\; \mathfrak{M} = 2+\sqrt{2},
\end{array}
$$

\noindent
because $t = \displaystyle\frac{p}{q - p} > 0$ for $\mbox{\small $\displaystyle\frac{1}{2t}$} + t  \geq \sqrt{2}$.
\hfill $\Box$ \end{proof}

\begin{proposition}
Let $D{\big (}p,\frac{r}{q}(q-p){\big )} \in \left[\Pi_{12}\right]$. If the inequality $(\ref{E8})$ holds for $D{\big (}p,\frac{r}{q}(q-p){\big )} \in \left[\Pi_{12}\right]$, then the following conclusion holds for the weight coefficient $w\!:$
  \begin{equation}
\langle \textbf{\textit{b}} \rangle  \; 0 \!\leq\! p \!<\! r\! \leq\! q
\;\Longrightarrow\;
      w \leq \omega(p,q,r) = 1 + \displaystyle\frac{q-p}{r+p} + \displaystyle\frac{q}{q-p}
  \end{equation}
and in that case
  \begin{equation}
  \omega(p,q,r) \in [\mathfrak{M}, \infty) \;\;\mbox{and}\;\; \mathfrak{M}=\frac{3}{2}+\sqrt{2}.
  \end{equation}
 \end{proposition}
\begin{proof}
Let us consider  $\langle \textbf{\textit{b}} \rangle  \; 0 \!\leq\! p \!<\! r\! \leq\! q$. Then, we notice the following expression holds:
$$
\begin{array}{rclc}
2\,\omega(p,q,r)
\!&\!=\!&\!
2 + 2\,\displaystyle\frac{q-p}{r+p} + \displaystyle\frac{2\,q}{q-p}            &                                      \\[2.5 ex]
\!&\!=\!&\!
3 + 2\,\displaystyle\frac{q-p}{r+p} + \displaystyle\frac{q+p}{q-p}             &                                      \\[2.5 ex]
\!&\!\geq\!&\!
3 + 2\,\displaystyle\frac{q-p}{q+p} + \displaystyle\frac{q+p}{q-p}  \geq 3 + 2\sqrt{2} & \;\Longrightarrow\; \mathfrak{M} = \mbox{\small $\displaystyle\frac{3}{2}$}+\sqrt{2},
\end{array}
$$
because $t = \displaystyle\frac{q-p}{q+p} > 0$ holds $2\,t + \mbox{\small $\displaystyle\frac{1}{t}$} \geq 2\sqrt{2}$.
\hfill $\Box$ \end{proof}

\begin{proposition}
Let $D{\big (}p,\frac{r}{q}(q-p){\big )} \in \left[\Pi_{12}\right]$. If the inequality $(\ref{E8})$ holds for $D{\big (}p,\frac{r}{q}(q-p){\big )} \in \left[\Pi_{12}\right]$, then the following conclusion holds for the weight coefficient $w\!:$
  \begin{equation}
\langle \textbf{\textit{c}} \rangle  \; 0 \!\leq\! p \!<\! q\! < \! r
\;\Longrightarrow\;
      w \leq \omega(p,q,r) = 1 + \displaystyle\frac{q-p}{r+p} + \displaystyle\frac{q}{q-p}
  \end{equation}
and in that case
  \begin{equation}
  \omega(p,q,r) \in (\mathfrak{M}, \infty) \;\;\mbox{and}\;\; \mathfrak{M}=2.
  \end{equation}
 \end{proposition}
\begin{proof}
Let us consider $\langle \textbf{\textit{c}} \rangle  \; 0 \!\leq\! p \!<\! q\! < \! r$. Then, we notice the following expression holds:
$$
\omega(p,q,r)
=
1 + \displaystyle\frac{q-p}{r+p} + \displaystyle\frac{q}{q-p}
\geq
2 + \displaystyle\frac{q-p}{r+p}
>
2
\;\Longrightarrow\; \mathfrak{M} = 2,
$$
because $\mbox{\small $\displaystyle\frac{q}{q-p}$} \geq 1$.
\hfill $\Box$ \end{proof}

\bigskip


Similar to previous propositions, the following three propositions are obtained from Statement 8.

\begin{proposition}
Let $O(0,0)\in \left[\Pi_{22}\right]$. If the inequality $(\ref{E8})$ holds for $O(0,0)$, then the following conclusion holds for the weight coefficient $w\!:$
  \begin{equation}
\langle \textbf{\textit{a}} \rangle  \; 0 \!<\! r \!\leq\! -p\! \leq \! q
\;\Longrightarrow\;
      w \leq \omega(p,q,r) = \frac{1}{2} + \frac{q-p}{2r}
  \end{equation}
and in that case
  \begin{equation}
  \omega(p,q,r) \in [\mathfrak{M}, \infty) \;\;\mbox{and}\;\; \mathfrak{M}=\mbox{\small $\displaystyle\frac{3}{2}$}.
  \end{equation}
 \end{proposition}
\begin{proof}
Let us consider $\langle \textbf{\textit{a}} \rangle  \; 0 \!<\! r \!\leq\! -p\! \leq \! q $. Then, we notice the following expression holds:
$$
\omega(p,q,r)
=
\frac{1}{2} + \frac{q}{2r} + \frac{-p}{2r}
\geq
\frac{1}{2} + \frac{q}{2(-p)} + \frac{-p}{2(-p)}
\geq
\frac{3}{2}
\;\Longrightarrow\;\; \mathfrak{M}=\frac{3}{2}.
$$
\hfill $\Box$ \end{proof}

\begin{proposition}
Let $O(0,0)\in \left[\Pi_{22}\right]$. If the inequality $(\ref{E8})$ holds for $O(0,0)$, then the following conclusion holds for the weight coefficient $w\!:$
  \begin{equation}
\langle \textbf{\textit{b}} \rangle  \; 0 \! < \! -p \! \leq \! r \! \leq \! q
\;\Longrightarrow\;
      w \leq \omega(p,q,r) = 1 + \displaystyle\frac{q}{r-p}
  \end{equation}
and in that case
  \begin{equation}
  \omega(p,q,r) \in [\mathfrak{M}, \infty) \;\;\mbox{and}\;\; \mathfrak{M}=\mbox{\small $\displaystyle\frac{3}{2}$}.
  \end{equation}
 \end{proposition}
\begin{proof}
Let us consider $\langle \textbf{\textit{b}} \rangle  \; 0 \! < \! -p \! \leq \! r \! \leq \! q$. Then, we notice the following expression holds:
$$
\omega(p,q,r)
=
1 + \displaystyle\frac{q}{r-p}
\geq
1 + \displaystyle\frac{q}{2r}
\geq
1 + \displaystyle\frac{1}{2} =  \frac{3}{2} \;\;\Longrightarrow\;\; \mathfrak{M} = \frac{3}{2}.
$$
\hfill $\Box$ \end{proof}

\begin{proposition}
Let $O(0,0)\in \left[\Pi_{22}\right]$. If the inequality $(\ref{E8})$ holds for $O(0,0)$, then the following conclusion holds for the weight coefficient $w\!:$
  \begin{equation}
\langle \textbf{\textit{c}} \rangle  \; 0 \! < \! -p \! \leq \! q \!<\! r
\;\Longrightarrow\;
      w \leq \omega(p,q,r) = 1 + \displaystyle\frac{r}{q-p}
  \end{equation}
and in that case
  \begin{equation}
  \omega(p,q,r) \in (\mathfrak{M}, \infty) \;\;\mbox{and}\;\; \mathfrak{M}=\mbox{\small $\displaystyle\frac{3}{2}$}.
  \end{equation}
 \end{proposition}
\begin{proof}
Let us consider $\langle \textbf{\textit{c}} \rangle  \; 0 \! < \! -p \! \leq \! q \!<\! r$. Then, we notice the following expression holds:
$$
\omega(p,q,r)
=
1 + \displaystyle\frac{r}{q-p}
>
1 + \displaystyle\frac{q}{q-p}
\geq
1 + \displaystyle\frac{q}{2q}
=\frac{3}{2} \;\;\Longrightarrow\;\; \mathfrak{M} = \frac{3}{2}.
$$
\hfill $\Box$ \end{proof}

Let us emphasize that the results of the previous three Propositions provide an improvement over some results from paper \cite{Ghandehari_2018}.


\section{Summa summarum}

Based on the propositions above, a theorem follows:
\begin{theorem}
In taxicab geometry for an interior point of $\triangle ABC$ in an appropriate position, the {\sc Erd\" os-Mordell}'s inequality holds
$$
R_A + R_B + R_C
\geq
\frac{3}{2}\left(r_a + r_b + r_c\right)_{\mbox{\micros \symbol{5}}}
$$
\end{theorem}

It is well known that taxicab distance depends on the rotation of the coordinate system, but does not depend on its translation or its reflection over a coordinate axis \cite{Ozcan_Ekmekci_Bayar_2002}. For an arbitrary triangle $ABC$ we set the following open problem (illustrated by Fig. 5).

\begin{conjecture}
In taxicab geometry for an interior point of any triangle $ABC$ the {\sc Erd\" os-Mordell}'s inequality holds
$$
R_A + R_B + R_C
\geq
\frac{3}{2}\left(r_a + r_b + r_c\right)_{\mbox{\micros \symbol{5}}}
$$
\end{conjecture}

\smallskip
\begin{center} 

\includegraphics*[height=55.0mm]{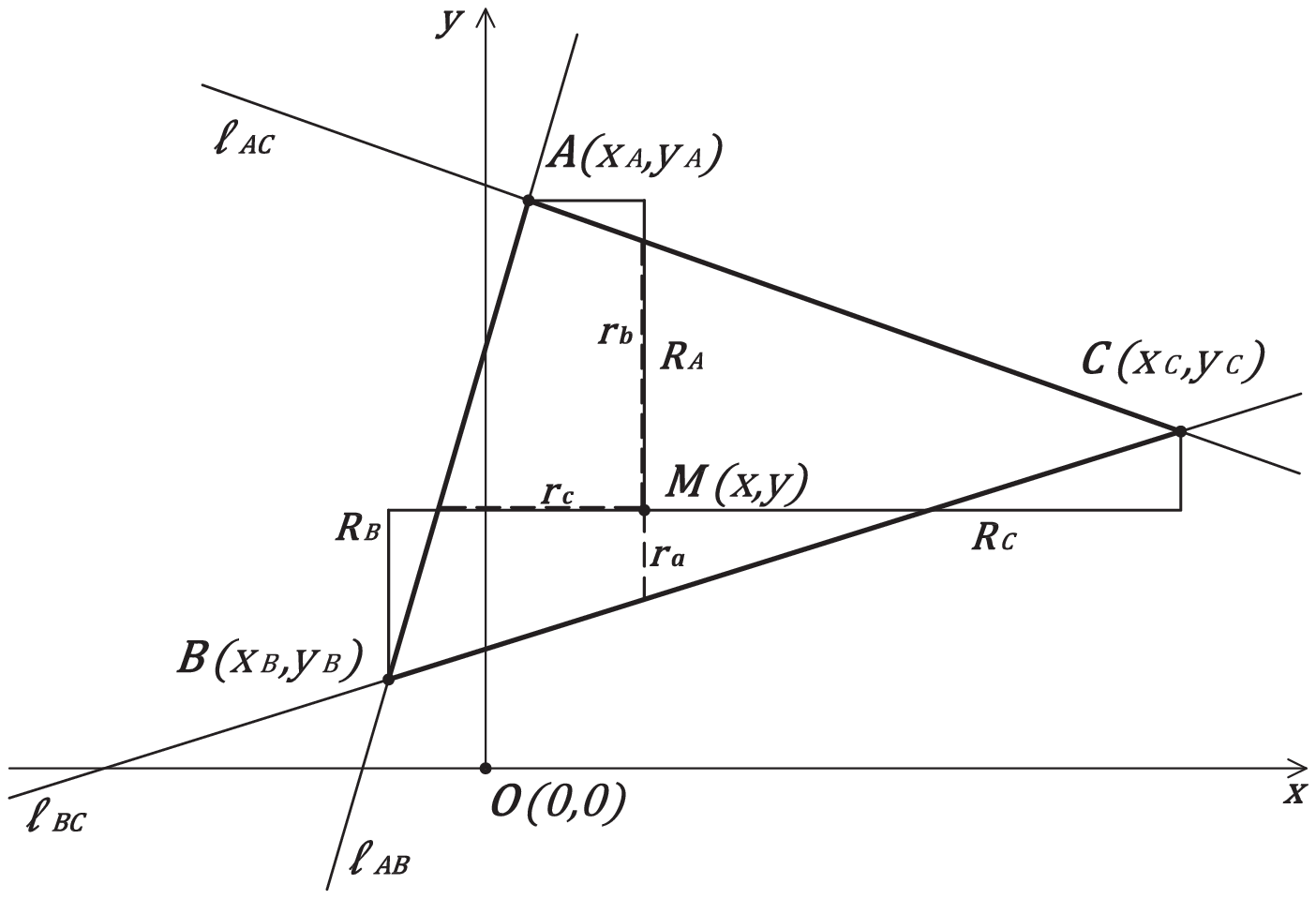}

\smallskip
\noindent
\textit{Figure 5. A geometric illustration of conjecture 1}
\end{center}

\bigskip

\bigskip
{\bf Acknowledgment.} Research of the second author was supported in part by the Serbian Ministry of Education, Science and Technological Development, under Projects ON 174032 and III 44006.

\end{document}